\documentclass{birkjour}
 \newtheorem{thm}{Theorem}[section]
 \newtheorem{cor}[thm]{Corollary}
 \newtheorem{lem}[thm]{Lemma}
 \newtheorem{prop}[thm]{Proposition}
 \newtheorem{conj}[thm]{Conjecture}
 \theoremstyle{definition}
 \newtheorem{defn}[thm]{Definition}
 \theoremstyle{remark}
 \newtheorem{rem}[thm]{Remark}
 \newtheorem*{ex}{Example}
 \numberwithin{equation}{section}
 
 \usepackage{mathrsfs}
 \usepackage{bm}
 \usepackage{calligra} 
 \usepackage[all,cmtip]{xy}
 \usepackage[OT2,T1]{fontenc}
 \usepackage{mathtools}
 
 \newcommand{\beq}{\begin{equation}}
 \newcommand{\eneq}{\end{equation}}
 \newcommand{\Om}{\Omega}
 \newcommand{\Hom}{{\rm Hom}}

 \newlanguage\fakelanguage
\newcommand\cyr{\fontencoding{OT2}\fontfamily{wncyr}\selectfont
   \language\fakelanguage}
\DeclareTextFontCommand{\textcyr}{\cyr}

\begin{document}

%-------------------------------------------------------------------------
% editorial commands: to be inserted by the editorial office
%
%\firstpage{1} \volume{228} \Copyrightyear{2004} \DOI{003-0001}
%\seriesextraline{This is the Concrete Title of this Book\br H.E. R and S.T.C. W, Eds.}
%\firstpage{1}

\title[Quantum differential equations and helices]
 {Quantum differential equations and helices}

%----------Author 1
\author[Giordano Cotti]{Giordano Cotti}

\address{%
Max-Planck Institut f\"{u}r Mathematik\\
Vivatsgasse 7\\
53111 Bonn\\
Germany}

\email{gcotti@sissa.it, gcotti@mpim-bonn.mpg.de}

%\thanks{This work was completed with the support of our
%----------Author 2
%\author{A Second Author}
%\address{The address of\br
%the second author\br
%sitting somewhere\br
%in the world}
%\email{dont@know.who.knows}
\subjclass{%Primary 
53D45; %Secondary 
18E30}

\keywords{Quantum cohomology, Frobenius manifolds, monodromy data, exceptional collections, Dubrovin's conjecture.}

\date{January 1, 2004}
%\dedicatory{To Boris Dubrovin, with admiration and gratitude}
%%% ----------------------------------------------------------------------

\begin{abstract}
These notes are a short and self-contained introduction to the isomonodromic approach to quantum cohomology, and Dubrovin's conjecture. An overview of recent results obtained in joint works with B.\,Dubrovin and D.\,Guzzetti \cite{CDG2}, and A.\,Varchenko \cite{CV} is given.
\end{abstract}

%%% ----------------------------------------------------------------------
\maketitle
%%% ----------------------------------------------------------------------
%\tableofcontents
\section{Quantum cohomology}
\subsection{Notations and conventions} Let $X$ be a smooth projective variety over $\mathbb C$ with vanishing odd-cohomology, i.e. $H^{2k+1}(X,\mathbb C)=0$, for $k\geq 0$. Fix a homogeneous basis $(T_1,\dots, T_n)$ of the complex vector space $H^\bullet(X):=\bigoplus_k H^{2k}(X,\mathbb C)$, and denote by $\bm t:=(t^1,\dots, t^n)$ the corresponding dual coordinates. Without loss of generality, we assume that $T_1=1$. The \emph{Poincar\'e pairing} on $H^\bullet(X)$ will be denoted by
\beq
\eta(u,v):=\int_Xu\cup v,\quad u,v\in H^\bullet(X),
\eneq
and we put $\eta_{\alpha\beta}:=\eta(T_\alpha, T_\beta)$, for $\alpha,\beta=1,\dots, n$, to be the Gram matrix wrt the fixed basis. The entries of the inverse matrix will be denoted by $\eta^{\alpha\beta}$, for $\alpha,\beta=1,\dots, n$.
In all the paper, the Einstein rule of summation over repeated indices is used. General references for this Section are \cite{CDG1,CDG2,Dnap,D1,D0,D2,KM,M,RT}.

\subsection{Gromov-Witten invariants in genus 0}\label{gwsec}
For a fixed $\beta\in H_2(X,\mathbb Z)/{\rm torsion}$, denote by $\overline{\mathcal M}_{0,k}(X,\beta)$ the Deligne-Mumford moduli stack of $k$-pointed stable rational maps with target $X$ of degree $\beta$:
\beq
\overline{\mathcal M}_{0,k}(X,\beta):=\left\{f\colon (C,\bm x)\to X,\ f_*[C]=\beta\right\}/{\rm equivalencies,}
\eneq
where $C$ is an algebraic curve of genus 0 with at most nodal singularities, $\bm x:=(x_1,\dots, x_k)$ is a $k$-tuple of pairwise distinct marked points of $C$, and equivalencies are automorphisms of $C\to X$ identical on $X$ and the markings. 

\emph{Gromov-Witten invariants} ($GW$-invariants for short)  of $X$, and their \emph{descendants}, are  defined as intersection numbers of cycles on $\overline{\mathcal M}_{0,k}(X,\beta)$, by the integrals
\beq\label{gw}
\langle\tau_{d_1}\gamma_1,\dots,\tau_{d_k}\gamma_k\rangle_{k,\beta}^X:=\int_{[\overline{\mathcal M}_{0,k}(X,\beta)]^{\rm virt}}\prod_{i=1}^k{\rm ev}_i^*\gamma_i\wedge \psi_i^{d_i},
\eneq
for $\gamma_1,\dots,\gamma_k\in H^\bullet(X)$, $d_i\in\mathbb N$. In formula \eqref{gw}, \beq{\rm ev}_i\colon \overline{\mathcal M}_{0,k}(X,\beta)\to X,\quad f\mapsto f(x_i),\quad i=1,\dots,k,
\eneq are evaluation maps, and $\psi_i:=c_1(\mathcal L_i)$ are the first Chern classes of the universal cotangent line bundles 
\beq
\mathcal L_i\to \overline{\mathcal M}_{0,k}(X,\beta),\quad \mathcal L_i|_{f}=T_{x_i}^*C,\quad i=1,\dots,k.
\eneq
The \emph{virtual fundamental cycle} $[\overline{\mathcal M}_{0,k}(X,\beta)]^{\rm virt}$ is an element of the Chow ring $A_\bullet\left(\overline{\mathcal M}_{0,k}(X,\beta)\right)$, namely
\[
[\overline{\mathcal M}_{0,k}(X,\beta)]^{\rm virt}\in A_D\left(\overline{\mathcal M}_{0,k}(X,\beta)\right),\quad D:=\dim_{\mathbb C}X-3+k+\int_\beta c_1(X).
\]
See \cite{BF} for its construction.
\subsection{Quantum cohomology as a Frobenius manifold}

Introduce infinitely many variables $\bm t_\bullet:=(t^\alpha_p)_{\alpha,p}$ with $\alpha=1,\dots, n$ and $p\in\mathbb N$. 

\begin{defn}
The \emph{genus 0 total descendant potential} of $X$ is the generating function $\mathcal F_0^X\in\mathbb C[\![\bm t_\bullet]\!]$ of descendant $GW$-invariants of $X$ defined by
\[
\mathcal F_0^X(\bm t_\bullet):=\sum_{k=0}^\infty\sum_{\beta}\sum_{\alpha_1,\dots,\alpha_k=1}^n\sum_{p_1,\dots, p_k=0}^\infty\frac{t^{\alpha_1}_{p_1}\dots t^{\alpha_k}_{p_k}}{k!}\langle \tau_{p_1}T_{\alpha_1},\dots, \tau_{p_k}T_{\alpha_k}\rangle_{k,\beta}^X.
\]
Setting $t^\alpha_0=t^\alpha$ and $t^\alpha_p=0$ for $p>0$, we obtain the \emph{Gromov-Witten potential} of $X$
\beq\label{gwpot}
F^X_0(\bm t):=\sum_{k=0}^\infty\sum_{\beta}\sum_{\alpha_1,\dots,\alpha_k=1}^n\frac{t^{\alpha_1}\dots t^{\alpha_k}}{k!}\langle T_{\alpha_1},\dots, T_{\alpha_k}\rangle_{k,\beta}^X.
\eneq
\end{defn}

Let $\Omega\subseteq H^\bullet(X)$ be the domain of convergence of $F_0^X(\bm t)$, assumed to be non-empty. 
We denote by $T\Omega$ and $T^*\Omega$ its holomorphic tangent and cotangent bundles, respectively.
Each tangent space $T_p\Omega$, with $p\in\Omega$, is canonically identified with the space $H^\bullet(X)$, via the identification $\frac{\partial}{\partial t^\alpha}\mapsto T_\alpha$.
The Poincar\'e metric $\eta$ defines a flat non-degenerate $\mathcal O_\Omega$-bilinear pseudo-riemannian metric on $\Omega$. The coordinates $\bm t$ are manifestly flat.  
Denote by $\nabla$ the Levi-Civita connection of $\eta$.

\begin{defn}Define the tensor $c\in \Gamma(T\Omega\otimes \bigodot^2T^*\Omega)$ by
\beq
c^\alpha_{\beta\gamma}:=\eta^{\alpha\lambda}\nabla^3_{\lambda\beta\gamma}F^X_0,\quad\alpha,\beta,\gamma=1,\dots, n,
\eneq
and let us introduce a product $*$ on vector fields on $\Omega$ by
\beq
\frac{\partial}{\partial t^\beta}*\frac{\partial}{\partial t^\gamma}:=c^\alpha_{\beta\gamma}\frac{\partial}{\partial t^\alpha},\quad \beta,\gamma=1,\dots, n.
\eneq
\end{defn}

\begin{thm}[\cite{KM,RT}]
The Gromov-Witten potential $F^X_0(\bm t)$ is a solution of $WDVV$ equations
\beq
\frac{\partial^3 F^X_0(\bm t)}{\partial t^\alpha\partial t^\beta\partial t^\gamma}\eta^{\gamma\delta}\frac{\partial^3 F^X_0(\bm t)}{\partial t^\delta\partial t^\epsilon\partial t^\phi}=\frac{\partial^3 F^X_0(\bm t)}{\partial t^\phi\partial t^\beta\partial t^\gamma}\eta^{\gamma\delta}\frac{\partial^3 F^X_0(\bm t)}{\partial t^\delta\partial t^\epsilon\partial t^\alpha},
\eneq
for $\alpha,\beta,\epsilon,\phi=1,\dots, n$.
\end{thm}

On each tangent space $T_p\Omega$, the product $*_p$ defines a structure of associative, commutative algebra with unit $\frac{\partial}{\partial t^1}\equiv 1$. Furthermore, the product $*$ is compatible with the Poincar\'e metric, namely
\beq
\eta(u*v,w)=\eta(u,v*w),\quad u,v,w\in\Gamma(T\Omega).
\eneq
This endows $(T_p\Omega, *_p,\eta_p,\left.\frac{\partial}{\partial t^1}\right|_p)$ with a complex \emph{Frobenius algebra} structure.

\begin{defn}
The vector field
\beq
E=c_1(X)+\sum_{\alpha=1}^n\left(1-\frac{1}{2}\deg T_\alpha\right)t^\alpha\frac{\partial}{\partial t^\alpha},
\eneq
is called  \emph{Euler vector field}. Here, $\deg T_\alpha$ denotes the cohomological degree of $T_\alpha$, i.e. $\deg T_\alpha:=r_\alpha$ if and only if $T_\alpha\in H^{r_\alpha}(X,\mathbb C).$
We denote by $\mathcal U$ the $(1,1)$-tensor defined by the multiplication with the Euler vector field, i.e.
\beq
\mathcal U\colon \Gamma(T\Omega)\to \Gamma(T\Omega),\quad v\mapsto E*v.
\eneq
\end{defn}

\begin{prop}[\cite{D1,D2}]
The Euler vector field $E$ is a Killing conformal vector field, whose flow preserves the structure constants of the Frobenius algerbas:
\beq
\mathfrak L_E\eta=(2-\dim_\mathbb CX)\eta,\quad \mathfrak L_Ec=c.
\eneq
\end{prop}

The structure $(\Om,c,\eta,\frac{\partial}{\partial t^1},E)$ gives an example of analytic \emph{Frobenius manifold}, called \emph{quantum cohomology of $X$} and denoted by $QH^\bullet(X)$, see \cite{D1,D0,D2,M}.
\subsection{Extended deformed connection}

\begin{defn}
The \emph{grading operator} $\mu\in {\rm End}(T\Omega)$ is the tensor  defined by
\beq
\mu(v):=\frac{2-\dim_\mathbb CX}{2}v-\nabla_vE,\quad v\in\Gamma(T\Omega).
\eneq
\end{defn}

Consider the canonical projection $\pi\colon\mathbb C^*\times \Omega\to \Omega$, and  the pull-back bundle $\pi^*T\Omega$. Denote by
\begin{enumerate}
\item $\mathscr T_\Om$ the sheaf of sections of $T\Om$, 
\item $\pi^*\mathscr T_\Om$ the pull-back sheaf, i.e. the sheaf of sections of $\pi^*T\Om$ 
\item  $\pi^{-1}\mathscr T_\Om$ the sheaf of sections of $\pi^*T\Om$ constant on the fibers of $\pi$.
\end{enumerate} 
All the tensors $\eta,c, E, \mathcal U,\mu$ can be lifted to $\pi^*T\Omega$, and their lifts will be denoted by the same symbols. The Levi-Civita connection $\nabla$ is lifted on $\pi^*T\Om$, and it acts so that
\beq
\nabla_\frac{\partial}{\partial z}v=0\quad \text{for }v\in(\pi^{-1}\mathscr T_\Om)(\Om),
\eneq
where $z$ is the coordinate on $\mathbb C^*$.

\begin{defn}
The \emph{extended deformed connection} is the connection $\widehat \nabla$ on the bundle $\pi^*T\Omega$ defined by
\begin{align}
\widehat\nabla_{w}v&=\nabla_wv+z\cdot w*v,\\
\widehat\nabla_{\frac{\partial}{\partial z}}v&=\nabla_{\partial_z}v+\mathcal U(v)-\frac{1}{z}\mu(v),
\end{align}
for $v,w\in \Gamma(\pi^*T\Omega)$.
\end{defn}

\begin{thm}[\cite{D1,D2}]
The connection $\widehat\nabla$ is flat.
\end{thm}

\subsection{Semisimple points and orthonormalized idempotent frame} 
\begin{defn}
A point $p\in \Om$ is \emph{semisimple} if and only if the corresponding Frobenius algebra $(T_p\Om,*_p,\eta_p, \frac{\partial}{\partial t^1}|_p)$ is without nilpotents. Denote by $\Om_{ss}$ the open dense subset of $\Om$ of semisimple points.
\end{defn}

\begin{thm}[\cite{HMT}]
The set $\Om_{ss}$ is non-empty only if $X$ is of Hodge-Tate\footnote{Here $h^{p,q}(X):=\dim_{\mathbb C}H^q(X,\Om_X^p)$, with $\Om_X^p$ the sheaf of holomorphic $p$-forms on $X$, denotes the $(p,q)$-Hodge number of $X$.} type, i.e. $h^{p,q}(X)=0$ for $p\neq q$.
\end{thm}
On $\Om_{ss}$ there are $n$ well-defined idempotent vector fields $\pi_1,\dots,\pi_n\in\Gamma(T\Om_{ss})$, satisfying
\beq
\pi_i*\pi_j=\delta_{ij}\pi_i,\quad\eta(\pi_i,\pi_j)=\delta_{ij}\eta(\pi_i,\pi_i),\quad i,j=1,\dots,n.
\eneq

\begin{thm}[\cite{Dnap,D1,D2}]\label{canu}
The idempotent vector fields pairwise commute: $[\pi_i,\pi_j]=0$ for $i,j=1,\dots, n$. Hence, there exist holomorphic local coordinates $(u_1,\dots, u_n)$ on $\Om_{ss}$ such that $\frac{\partial}{\partial u_i}=\pi_i$ for $i=1,\dots, n$. 
\end{thm}

\begin{defn}
The coordinates $(u_1,\dots, u_n)$ of Theorem \ref{canu} are called \emph{canonical coordinates}. 
\end{defn}

\begin{prop}[\cite{D1,D2}]
Canonical coordinates are uniquely defined up to ordering and shifts by constants. The eigenvalues of the tensor $\mathcal U$ define a system of canonical coordinates in a neighborhood of any semisimple point of $\Om_{ss}$.
\end{prop}

\begin{defn}
We call \emph{orthonormalized idempotent frame} a frame $(f_i)_{i=1}^n$ of $T\Om_{ss}$  defined by
\beq\label{fvects}
f_i:=\eta(\pi_i,\pi_i)^{-\frac{1}{2}}\pi_i,\quad i=1,\dots,n,
\eneq 
for arbitrary choices of signs of the square roots. The $\Psi$-matrix is the matrix $(\Psi_{i\alpha})_{i,\alpha=1}^n$ of change of tangent frames,  defined by
 \beq
 \frac{\partial}{\partial t^\alpha}=\sum_{i=1}^n\Psi_{i\alpha}f_i,\quad \alpha=1,\dots,n.
 \eneq
\end{defn}

\begin{rem}
In the orthonormalized idempotent frame, the operator $\mathcal U$ is represented by a diagonal matrix, and the operator $\mu$ by an antisymmetric matrix:
\beq
U:=\operatorname{diag}(u_1,\dots,u_n),\quad\Psi\mathcal U\Psi^{-1}=U,
\eneq
\beq
V:=\Psi\mu\Psi^{-1},\quad V^T+V=0.
\eneq
\end{rem}

\section{Quantum differential equation} The connection $\widehat\nabla$ induces a flat connection on $\pi^*(T^*\Om)$. Let $\xi\in\Gamma(\pi^*(T^*\Om))$ be a flat section. Consider the corresponding vector field $\zeta\in\Gamma(\pi^*(T\Om))$ via musical isomorphism, i.e. such that $\xi(v)=\eta(\zeta,v)$ for all $v\in\Gamma(\pi^*(T\Om))$.

The vector field $\zeta$ satisfies the following system\footnote{We consider the joint system \eqref{eq1}, \eqref{qde} in matrix notations ($\zeta$ a column vector whose entries are the components $\zeta^\alpha(\bm t,z)$ wrt $\frac{\partial}{\partial t^\alpha}$). Bases of solutions are arranged in invertible $n\times n$-matrices, called \emph{fundamental systems of solutions}.} of equations
 \begin{align}
 \label{eq1}
 \frac{\partial}{\partial t^\alpha}\zeta&=z\mathcal C_\alpha\zeta,\quad \alpha=1,\dots, n,\\
\label{qde}
  \frac{\partial}{\partial z}\zeta&=\left(\mathcal U+\frac{1}{z}\mu\right)\zeta.
 \end{align}
Here $\mathcal C_\alpha$ is the $(1,1)$-tensor defined by $(\mathcal C_\alpha)^\beta_\gamma:=c^\beta_{\alpha\gamma}$. 
\begin{defn}
The \emph{quantum differential equation} ($qDE$) of $X$ is the differential equation \eqref{qde}. 
\end{defn}

The $qDE$ is an ordinary differential equation with rational coefficients. It has two singularities on the Riemann sphere $\mathbb P^1(\mathbb C)$:
\begin{enumerate}
\item a Fuchsian singularity at $z=0$,
\item an irregular singularity (of Poincar\'{e} rank 1) at $z=\infty$.
\end{enumerate}
Points of $\Om$ are parameters of deformation of the coefficients of the $qDE$. Solutions $\zeta(\bm t,z)$ of the joint system of equations \eqref{eq1}, \eqref{qde} are ``multivalued'' functions wrt $z$, i.e. they are well-defined functions on $\Om\times\widehat{\mathbb C^*}$, where $\widehat{\mathbb C^*}$ is the universal cover of $\mathbb C^*$.

\subsection{Solutions in Levelt form at $z=0$ and topological-enumerative solution}
\begin{thm}[\cite{CDG1,D1,D2}]\label{levelt}
There exist fundamental systems of solutions $Z_0(\bm t,z)$ of the joint system \eqref{eq1}, \eqref{qde} with expansions at $z=0$ of the form
\beq\label{levelt}
Z_0(\bm t,z)=F(\bm t,z)z^\mu z^R,\quad R=\sum_{k\geq 1} R_k, \quad F(\bm t,z)=I+\sum_{j=1}^\infty F_j(\bm t)z^j
\eneq
where $(R_{k})_{\alpha\beta}\neq 0\text{ only if }\mu_\alpha-\mu_\beta=k$. The series $F(\bm t,z)$ is convergent and  satisfies the orthogonality condition
\beq
F(\bm t,-z)^T\eta F(\bm t,z)=\eta.
\eneq
\end{thm}

\begin{defn}
A fundamental system of solutions $Z_0(\bm t,z)$ of the form described in Theorem \ref{levelt} are said to be in \emph{Levelt form} at $z=0$.
\end{defn}

\begin{rem}\label{Rc1}
Fundamental systems of solutions in Levelt form are not unique. The exponent $R$ is not uniquely determined. Moreover, even for a fixed exponent $R$, the series $F(\bm t,z)$ is not uniquely determined, see \cite{CDG1}. It can be proved that the matrix $R$ can be chosen as the matrix of the operator $c_1(X)\cup (-)\colon H^\bullet(X)\to H^\bullet(X)$ wrt the basis $(T_\alpha)_{\alpha=1}^n$ \cite[Corollary 2.1]{D2}.
\end{rem}

\begin{rem}\label{M0}
Let $Z_0(\bm t, z)$ be a fundamental system of solutions in Levelt form \eqref{levelt}. The monodromy matrix $M_0(\bm t)$, defined by
\beq
Z_0(\bm t, e^{2\pi\sqrt{-1}}z)=Z_0(\bm t,z)M_0(\bm t),\quad z\in\widehat{\mathbb C^*},
\eneq is given by
\beq
M_0(\bm t)=\exp(2\pi\sqrt{-1}\mu)\exp(2\pi\sqrt{-1}R).
\eneq
In particular, $M_0$ does not depend on $\bm t$.
\end{rem}

\begin{defn}Define the functions $\theta_{\beta,p}(\bm t, z),\,\theta_{\beta}(\bm t, z)$, with $\beta=1,\dots, n$ and $p\in\mathbb N$, by
\beq
\theta_{\beta,p}(\bm t):=\left.\frac{\partial^2\mathcal F_0^X(\bm t_\bullet)}{\partial t^1_0\partial t^\beta_p}\right|_{t^\alpha_p=0\text{ for }p>1,\quad t^\alpha_0=t^\alpha\text{ for }\alpha=1,\dots, n},
\eneq
\beq
\theta_\beta(\bm t,z):=\sum_{p=0}^\infty\theta_{\beta,p}(\bm t)z^p.
\eneq
Define the matrix $\Theta(\bm t,z)$ by
\beq
\Theta(\bm t,z)^\alpha_\beta:=\eta^{\alpha\lambda}\frac{\partial\theta_\beta(\bm t,z)}{\partial t^\lambda},\quad \alpha,\beta=1,\dots, n.
\eneq
\end{defn}

\begin{thm}[\cite{CDG1,D2}]
The matrix $Z_{\rm top}(\bm t,z):=\Theta(\bm t,z)z^\mu z^{c_1(X)\cup}$ is a fundamental system of solutions of the joint system \eqref{eq1}-\eqref{qde} in Levelt form at $z=0$.
\end{thm}

\begin{defn}
The solution $Z_{\rm top}(\bm t,z)$ is called \emph{topological-enumerative solution} of the joint system \eqref{eq1}, \eqref{qde}.
\end{defn}

\subsection{Stokes rays and $\ell$-chamber decomposition}
\begin{defn}
We call \emph{Stokes rays} at a point $p\in \Om$ the oriented rays $R_{ij}(p)$ in $\mathbb C$  defined by
\beq
R_{ij}(p):=\left\{-\sqrt{-1}(\overline{u_i(p)}-\overline{u_j(p)})\rho\colon \rho\in\mathbb R_{+}\right\},
\eneq
where $(u_1(p),\dots, u_n(p))$ is the spectrum of the operator $\mathcal U(p)$ (with a fixed arbitrary order).
\end{defn}
Fix an oriented ray $\ell$ in the universal cover $\widehat{\mathbb C^*}$.  
\begin{defn}
We say that $\ell$ is \emph{admissible} at $p\in\Om$ if the projection of the the ray $\ell$ on $\mathbb C^*$ does not coincide with any Stokes ray $R_{ij}(p)$.
\end{defn}
\begin{defn}
Define the open subset $O_\ell$ of points $p\in \Om$ by the following conditions:
\begin{enumerate}
\item the eigenvalues $u_i(p)$ are pairwise distinct,
\item $\ell$ is admissible at $p$.
\end{enumerate}
We call \emph{$\ell$-chamber} of $\Om$ any connected component of $O_\ell$. 
\end{defn}

\subsection{Stokes fundamental solutions at $z=\infty$}
Fix an oriented ray $\ell\equiv \left\{\arg z=\phi\right\}$ in $\widehat{\mathbb C^*}$. For $m\in\mathbb Z$, define the sectors in $\widehat{\mathbb C^*}$
\beq
\Pi_{L,m}(\phi):=\left\{z\in\widehat{\mathbb C^*}\colon \phi+2\pi m<\arg z<\phi+\pi+2\pi m\right\},
\eneq
\beq
\Pi_{R,m}(\phi):=\left\{z\in\widehat{\mathbb C^*}\colon \phi-\pi+2\pi m < \arg z< \phi+2\pi m\right\}.
\eneq
\begin{defn}The \emph{coalescence locus} of $\Om$ is the set
\beq
\Delta_\Om:=\left\{p\in\Om\colon u_i(p)=u_j(p),\quad \text{for some }i\neq j\right\}.
\eneq
\end{defn}
\begin{thm}[\cite{D1,D2}]
There exists a unique formal solution $Z_{\rm form}(\bm t,z)$ of the joint system \eqref{eq1}, \eqref{qde} of the form
\begin{align}
Z_{\rm form}(\bm t,z)&=\Psi(\bm t)^{-1}G(\bm t,z)\exp(z U(\bm t)),\\
G(\bm t,z)&= I+\sum_{k=1}^\infty\frac{1}{z^k}G_k(\bm t),
\end{align}
where the matrices $G_k(\bm t)$ are holomorphic on $\Om\setminus \Delta_\Om$.
\end{thm}
\begin{thm}[\cite{D1,D2}]
Let $m\in\mathbb Z$. There exist unique fundamental systems of solutions $Z_{L,m}(\bm t, z)$, $Z_{R,m}(\bm t, z)$ of the joint system \eqref{eq1}, \eqref{qde} with asymptotic expansion
\begin{align}
Z_{L,m}(\bm t,z)&\sim Z_{\rm form}(\bm t,z),\quad |z|\to\infty,\quad z\in \Pi_{L,m}(\phi),\\
Z_{R,m}(\bm t,z)&\sim Z_{\rm form}(\bm t,z),\quad |z|\to\infty,\quad z\in \Pi_{R,m}(\phi),
\end{align}
respectively.
\end{thm}

\begin{defn}
The solutions $Z_{L,m}(\bm t,z)$ and $Z_{R,m}(\bm t,z)$ are called \emph{Stokes fundamental solutions} of the joint system \eqref{eq1}, \eqref{qde} on the sectors $\Pi_{L,m}(\phi)$ and $\Pi_{R,m}(\phi)$ respectively.
\end{defn}

\subsection{Monodromy data}
Let $\ell\equiv \left\{\arg z=\phi\right\}$ be an oriented ray in $\widehat{\mathbb C^*}$ and consider the corresponding Stokes fundamental systems of solutions $Z_{L,m}(\bm t, z)$, $Z_{R,m}(\bm t, z)$, for $m\in\mathbb Z$.

\begin{defn}
We define the \emph{Stokes} and \emph{central connection} matrices $S^{(m)}(p)$, $C^{(m)}(p)$, with $m\in\mathbb Z$, at the point $p\in O_\ell$ by the identities 
\begin{align}
Z_{L,m}(\bm t(p), z)=Z_{R,m}(\bm t(p),z) S^{(m)}(p),\\
 Z_{R,m}(\bm t(p),z)=Z_{\rm top}(\bm t(p),z )C^{(m)}(p).
\end{align}
Set $S(p):=S^{(0)}(p)$ and $C(p):=C^{(0)}(p)$.
\end{defn}

\begin{defn}
The \emph{monodromy data} at the point $p\in O_\ell$ are defined as the $4$-tuple $(\mu,R,S(p),C(p))$, where
\begin{itemize}
\item $\mu$ is the (matrix associated to) the grading operator,
\item $R$ is the (matrix associated to) the operator $c_1(X)\cup\colon H^\bullet(X)\to H^\bullet(X)$,
\item $S(p), C(p)$ are the Stokes and central connection matrices at $p$, respectively.
\end{itemize}
\end{defn}

\begin{rem}
The definition of the Stokes and central connection matrices is subordinate to several non-canonical choices:
\begin{enumerate}
\item the choice of an oriented ray $\ell$ in $\widehat{\mathbb C^*}$,
\item the choice of an ordering of canonical coordinates $u_1,\dots, u_n$ on each $\ell$-chamber,
\item the choice of signs in \eqref{fvects}, and hence of the branch of the $\Psi$-matrix on each $\ell$-chamber.
\end{enumerate}
Different choices affect the numerical values of the data $(S,C)$, see \cite{CDG1}. In particular, for different choices of ordering of canonical coordinates, the Stokes and central connection matrices transform as follows:
\beq
S\mapsto \Pi S \Pi^{-1},\quad C\mapsto C\Pi^{-1},\quad \Pi\text{ permutation matrix}.
\eneq
\end{rem}

\begin{defn}
Fix a point $p\in O_\ell$ with canonical coordinates $(u_i(p))_{i=1}^n$. Define the oriented rays $L_j(p,\phi)$, $j=1,\dots,n$, in the complex plane by the equations
\beq
L_j(p,\phi):=\left\{u_j(p)+\rho e^{\sqrt{-1}(\frac{\pi}{2}-\phi)}\colon \rho\in\mathbb R_+\right\}.
\eneq
The ray $L_j(p,\phi)$ is oriented from $u_j(p)$ to $\infty$. We say that $(u_i(p))_{i=1}^n$ are in $\ell$-\emph{lexicographical order} if $L_j(p,\phi)$ is on the left of $L_k(p,\phi)$ for $1\leq j<k\leq n$. 
\end{defn}
In what follows, it is assumed that the $\ell$-lexicographical order of canonical coordinates is fixed at all points of $\ell$-chambers.
\begin{lem}[\cite{CDG1,D2}]
If the canonical coordinates $(u_i(p))_{i=1}^n$ are in $\ell$-lexicogra\-phical order at $p\in O_\ell$, then the Stokes matrices $S^{(m)}(p)$, $m\in\mathbb Z$, are upper triangular with $1$'s along the diagonal.
\end{lem}

By Remarks \ref{Rc1} and \ref{M0}, the matrices $\mu$ and $R$ determine the monodromy of solutions of the $qDE$,
\beq
M_0:=\exp(2\pi\sqrt{-1}\mu)\exp(2\pi\sqrt{-1}R).
\eneq Moreover, $\mu$ and $R$ do not depend on the point $p$. The following theorem furnishes a refinement of this property.

\begin{thm}[\cite{CDG1,D1,D2}]\label{thmd}
The monodromy data $(\mu,R, S,C)$ are constant in each $\ell$-chamber. Moreover, they satisfy the following identities:
\begin{align}
\label{const1}
CS^TS^{-1}C^{-1}&=M_0,\\
\label{const2}
S=C^{-1}\exp(-\pi\sqrt{-1}R)\exp&(-\pi\sqrt{-1}\mu)\eta^{-1}(C^T)^{-1},\\
\label{const3}
S^T=C^{-1}\exp(\pi\sqrt{-1}R)\exp&(\pi\sqrt{-1}\mu)\eta^{-1}(C^T)^{-1}.
\end{align}
\end{thm}

\begin{thm}[\cite{CDG1}]
The Stokes and central connection matrices $S_m,C_m$, with $m\in\mathbb Z$, can be reconstructed from the monodromy data $(\mu, R,S,C)$:
\beq\label{allsc}
S^{(m)}=S,\quad C^{(m)}=M_0^{-m}C,\quad m\in\mathbb Z.
\eneq
\end{thm}

\begin{rem}\label{coal}
Points of $O_\ell$ are semisimple. The results of \cite{CDG0,CDG1,CG1,CG2} imply that the monodromy data $(\mu,R,S,C)$ are well defined also at points $p\in\Om_{ss}\cap\Delta_{\Om}$, and that Theorem \ref{thmd} still holds true.
\end{rem} 
 \begin{rem}\label{RH}
 From the knowledge of the monodromy data $(\mu,R,S,C)$ the Gromov-Witten potential $F_0^X(\bm t)$ can be recostructed via a Riemann-Hilbert boundary value problem, see \cite{CDG1,CDG2,D2,G2}. Hence, the monodromy data may be interpreted as a \emph{system of coordinates} in the space of solutions of $WDVV$ equations.
 \end{rem}

\subsection{Action of the braid group $\mathcal B_n$}
Consider the braid group $\mathcal B_n$ with generators $\beta_1,\dots, \beta_{n-1}$ satisfying the relations 
\beq
\beta_{i}\beta_j=\beta_j\beta_i,\quad |i-j|>1,
\eneq
\beq
\beta_i\beta_{i+1}\beta_i=\beta_{i+1}\beta_i\beta_{i+1}.
\eneq
Let $\mathcal U_n$ be the set of upper triangular $(n\times n)$-matrices with $1$'s along the diagonal.
\begin{defn}Given $U\in\mathcal U_n$ define the matrices $A^{\beta_i}(U)$, with $i=1,\dots,n-1$, as follows
\begin{align}
\left(A^{\beta_i}(U)\right)_{hh}:=1,\quad h&=1,\dots,n,\quad h\neq i,i+1,\\
\left(A^{\beta_i}(U)\right)_{i+1,i+1}&=-U_{i,i+1},\\
\left(A^{\beta_i}(U)\right)_{i,i+1}&=\left(A^{\beta_i}(U)\right)_{i+1,i}=1,
\end{align}
and all other entries of $A^{\beta_i}(U)$ are equal to zero.
\end{defn}

\begin{lem}[\cite{CDG1,D1,D2}]\label{actbr}The braid group $\mathcal B_n$ acts on $\mathcal U_n\times GL(n,\mathbb C)$ as follows: 
\begin{align*}
\mathcal B_n\times \mathcal U_n\times GL(n,\mathbb C)&\xrightarrow{\hspace{2cm}}\mathcal U_n\times GL(n,\mathbb C)\\
(\beta_i,U,C)&\xmapsto{\quad\quad} (A^{\beta_i}(U)\cdot U\cdot A^{\beta_i}(U),\ C\cdot A^{\beta_i}(U)^{-1})
\end{align*}
We denote by $(U,C)^{\beta_i}$ the  action of $\beta_i$ on $(U,C)$.
\end{lem}

Fix an oriented ray $\ell\equiv\left\{\arg z=\phi\right\}$ in $\widehat{\mathbb C^*}$, and denote by $\overline{\ell}$ its projection on $\mathbb C^*$. Let $\Om_{\ell,1},\Om_{\ell,2}$ be two $\ell$-chambers and let $p_i\in\Om_{\ell,i}$ for $i=1,2$. The difference of values of the Stokes and central connection matrices $(S_1,C_1)$ and $(S_2,C_2)$, at $p_1$ and $p_2$ respectively, can be described by the action of the braid group $\mathcal B_n$ of Lemma \ref{actbr}. 

\begin{thm}[\cite{CDG1,D1,D2}]\label{thbr}
Consider a continuous path $\gamma\colon[0,1]\to\Om$ such that
\begin{itemize}
\item $\gamma(0)=p_1$ and $\gamma(1)=p_2$,
\item there exists a unique $t_o\in[0,1]$ such that $\ell$ is not admissible at $\gamma(t_o)$,
\item there exist $i_1,\dots, i_k\in\left\{1,\dots,n\right\}$, with $|i_a-i_b|>1$ for $a\neq b$, such that the rays\footnote{Here the labeling of Stokes rays is the one prolonged from the initial point $t=0$.} $\left(R_{i_j,i_j+1}(t)\right)_{j=1}^r$ (resp. $\left(R_{i_j,i_j+1}(t)\right)_{j=r+1}^k$) cross the ray $\overline{\ell}$ in the clockwise (resp. counterclockwise) direction, as $t\to t_o^-$.
\end{itemize}
Then, we have
\beq
(S_2,C_2)=(S_1,C_1)^\beta,\quad \beta:=\left(\prod_{j=1}^r\beta_{i_j}\right)\cdot\left(\prod_{h=r+1}^{k}\beta_{i_h}\right)^{-1}.
\eneq
\end{thm}
\begin{rem}
In the general case, the points $p_1$ and $p_2$ can be connected by concatenations of paths $\gamma$ satisfying the assumptions of Theorem \ref{thbr}. 
\end{rem}

\begin{rem}
The action of $\mathcal B_n$ on $(S,C)$ also describes the analytic continuation of the Frobenius manifold structure on $\Om$, see \cite[Lecture 4]{D2}.
\end{rem}

\section{Derived category, exceptional collections, helices}

\subsection{Notations and basic notions}
Denote by $Coh(X)$ the abelian category of coherent sheaves on $X$, and by $\mathcal D^b(X)$ its bounded derived category. Objects of $\mathcal D^b(X)$ are bounded complexes $A^\bullet$ of coherent sheaves on $X$. Morphisms are given by \emph{roofs}: if $A^\bullet, B^\bullet$ are two bounded complexes, a morphism $f\colon A^\bullet\to B^\bullet$ in $\mathcal D^b(X)$ is the datum of
\begin{itemize}
\item a third object $C^\bullet$ in $\mathcal D^b(X)$,
\item two homotopy classes of morphisms of complexes $q\colon C^\bullet \to A^\bullet$ and $g\colon C^\bullet \to B^\bullet$,
\item the morphism $q$ is required to be a \emph{quasi-isomorphism}, i.e. it induces isomorphism in cohomology.
\end{itemize} 
\beq
\xymatrix{&C^\bullet\ar[dl]_{q}\ar[dr]^{g}&\\
A^\bullet\ar@{-->}[rr]_{f}&&B^\bullet}
\eneq
The derived category $\mathcal D^b(X)$ admits a triangulated structure, the \emph{shift functor} $[1]\colon \mathcal D^b(X)\to \mathcal D^b(X)$ being defined by
\beq
A^\bullet[1]:=A^{\bullet+1},\quad A^\bullet\in\mathcal D^b(X).
\eneq
Denote by $\Hom^\bullet(A^\bullet, B^\bullet):=\bigoplus_{k\in\mathbb Z}\Hom(A^\bullet, B^\bullet[k])$. General references for this Section are \cite{GM,GK,GR,Rhel}.
\subsection{Exceptional collections}
\begin{defn}
An object $E\in \mathcal D^b(X)$ is called \emph{exceptional} iff 
\beq\Hom^\bullet(E,E)\cong \mathbb C.
\eneq
\end{defn}
\begin{defn}An \emph{exceptional  collection} is an ordered family $(E_1,\dots, E_n)$ of exceptional objects of $\mathcal D^b(X)$ such that
\beq \Hom^\bullet(E_j,E_i)\cong 0\quad\text{ for }j>i.
\eneq
An exceptional collection is \emph{full} if it generates $\mathcal D^b(X)$ as a triangulated category, i.e. if any full triangulated subcategory of $\mathcal D^b(X)$ containing all the objects $E_i$'s is equivalent to $\mathcal D^b(X)$ via the inclusion functor.
\end{defn}

\begin{ex}
In \cite{B} A.\,Beilinson showed that the collection of line bundles 
\beq\label{bec}\frak B:=(\mathcal O,\mathcal O(1),\dots,\mathcal O(n))\eneq 
on $\mathbb P^n$ is a full exceptional collection. M.\,Kapranov generalized this result in \cite{K}, where full exceptional collections on Grassmannians, flag varieties of group $SL_n$, and smooth quadrics are constructed. 

Denote by $\mathbb G(k,n)$ the Grassmannian of $k$-dimensional subspaces in $\mathbb C^n$, by $\mathcal S^\vee$ the dual of its tautological bundle. Let $\mathbb S^\lambda$ be the Schur functor (see \cite{F}) labelled by a Young diagram $\lambda$ inside a rectangle $k\times (n-k)$. The collection $\frak K:=\left(\mathbb S^\lambda\mathcal S^\vee\right)_{\lambda}$ is full and exceptional in $\mathcal D^b(\mathbb G(k,n))$. The order of the objects of the collection is the partial order defined by inclusion of Young diagrams.
\end{ex}

\subsection{Mutations and helices}
Let $E$ be an exceptional object in $\mathcal D^b(X)$. For any $X\in\mathcal D^b(X)$, we have natural evaluation and co-evaluation morphisms
\beq
j^*\colon \Hom^\bullet(E,X)\otimes E\to X,\quad j_*\colon X\to \Hom^\bullet(X,E)^*\otimes E.
\eneq
\begin{defn}
The \emph{left} and \emph{right mutations} of $X$ with respect to $E$ are the objects $\mathbb L_EX$ and $\mathbb R_EX$ uniquely defined by the distinguished triangles
\begin{equation}
\label{triangleleft}
\xymatrix{
\mathbb L_EX[-1]\ar[r]&\Hom^\bullet(E,X)\otimes E\ar[r]^{\quad\quad\quad j^*}&X\ar[r]&\mathbb L_EX,}
\end{equation}
\begin{equation}
\label{triangleright}
\xymatrix{
\mathbb R_EX\ar[r]&X\ar[r]^{ j_*\quad\quad\quad}&\Hom^\bullet(X,E)^*\otimes E\ar[r]&\mathbb R_EX[1],}
\end{equation}
respectively.
\end{defn}

\begin{rem}
In general, the third object of a distinguished triangle is not canonically defined by the other two terms. Nevertheless, the objects $\mathbb L_XE$ and $\mathbb R_EX$ are uniquely defined \emph{up to unique isomorphism}, because of the exceptionality of $E$, see \cite[Section 3.3]{CDG2}.
\end{rem}

\begin{defn}
Let $\frak E=(E_1,\dots, E_n)$ be an exceptional collection. For any $i=1,\dots, n-1$ define the \emph{left} and \emph{right mutations}
\begin{align}
\mathbb L_i\frak E:&=(E_1,\dots,\mathbb L_{E_i}E_{i+1},E_i,\dots, E_n),\\
\mathbb R_i\frak E:&=(E_1,\dots, E_{i+1},\mathbb R_{E_{i+1}}E_i,\dots, E_n).
\end{align}
\end{defn}

\begin{thm}[\cite{GK,Rhel}]\label{teobraid}
For all $i=1,\dots, n-1$ the collections $\mathbb L_i\frak E$ and $\mathbb R_i\frak E$ are exceptional. Moreover, we have that 
\begin{align*}
\mathbb L_i\mathbb R_i=\mathbb R_i\mathbb L_i={\rm Id},\quad &\mathbb L_{i+1}\mathbb L_i\mathbb L_{i+1}=\mathbb L_i\mathbb L_{i+1}\mathbb L_i,\quad i=1,\dots,n,\\
&\mathbb L_i\mathbb L_j=\mathbb L_j\mathbb L_i,\quad |i-j|>1.
\end{align*}
\end{thm}

According to Theorem \ref{teobraid}, we have a well-defined action of $\mathcal B_n$ on the set of exceptional collections of length $n$ in $\mathcal D^b(X)$: the action of the generator $\beta_i$ is identified with the action of the mutation $\mathbb L_i$ for $i=1,\dots, n-1$.

\begin{defn}
Let $\frak E=(E_1,\dots, E_n)$ be a full exceptional collection. We define the \emph{helix} generated by $\frak E$ to be the infinite family $(E_i)_{i\in\mathbb Z}$ of exceptional objects obtained by iterated mutations
\[
E_{n+i}:=\mathbb R_{E_{n+i-1}}\dots\mathbb R_{E_{i+1}}E_i,\quad E_{i-n}:=\mathbb L_{E_{i-n+1}}\dots\mathbb L_{E_{i-1}}E_i,\quad i\in\mathbb Z.
\]
Any family of $n$ consecutive exceptional objects $(E_{i+k})_{k=1}^n$ is called a \emph{foundation} of the helix.
\end{defn}

\begin{lem}[\cite{GK}]\label{sper}
For $i,j\in\mathbb Z$, we have $\Hom^\bullet(E_i,E_j)\cong \Hom^\bullet(E_{i-n},E_{j-n})$.
\end{lem}

\subsection{Exceptional bases in $K$-theory} Consider the Grothendieck group $K_0(X)\equiv K_0(\mathcal D^b(X))$, equipped with the Grothendieck-Euler-Poincar\'e bilinear form
\beq
\chi([V],[F]):=\sum_{k}(-1)^k\dim_{\mathbb C}\Hom(V,F[i]),\quad V,F\in\mathcal D^b(X).
\eneq
\begin{defn}
A basis $(e_i)_{i=1}^n$ of $K_0(X)_{\mathbb C}$ is called \emph{exceptional} if $\chi(e_i,e_i)=1$ for $i=1,\dots, n$, and $\chi(e_j,e_i)=0$ for $1\leq i<j\leq n$.
\end{defn}
\begin{lem}
Let $(E_i)_{i=1}^n$ be a full exceptional collection in $\mathcal D^b(X)$. The $K$-classes $([E_i])_{i=1}^n$ form an exceptional basis of $K_0(X)_{\mathbb C}$.
\end{lem}
The action of the braid group on the set of exceptional collections in $\mathcal D^b(X)$ admits a $K$-theoretical analogue on the set of exceptional bases of $K_0(X)_{\mathbb C}$, see \cite{CDG2,GK}.
\section{Dubrovin's conjecture} 
\subsection{$\Gamma$-classes and graded Chern character}Let $V$ be a complex vector bundle on $X$ of rank $r$, and let $\delta_1,\dots,\delta_r$ be its Chern roots, so that $c_j(V)=s_j(\delta_1,\dots, \delta_r)$, where $s_j$ is the $j$-th elementary symmetric polynomial.
\begin{defn}Let $Q$ be an indeterminate, and $F\in\mathbb C[\![Q]\!]$ be of the form $F(Q)=1+\sum_{n\geq 1}\alpha_nQ^n$. The $F$-\emph{class} of $V$ is the charcateristic class $\widehat F_V\in H^\bullet(X)$ defined by  $\widehat F_V:=\prod_{j=1}^rF(\delta_j).$
\end{defn}
\begin{defn}
The $\Gamma^{\pm}$-\emph{classes} of $V$ are the characteristic classes associated with the Taylor expansions
\beq
\Gamma(1\pm Q)=\exp\left(\mp\gamma Q+\sum_{m=2}^\infty(\mp1)^m\frac{\zeta(m)}{m}Q^n\right)\in\mathbb C[\![Q]\!],
\eneq
where $\gamma$ is the Euler-Mascheroni constant and $\zeta$ is the Riemann zeta function.
\end{defn}
If $V=TX$, then we denote $\widehat{\Gamma}_X^\pm$ its $\Gamma$-classes.
\begin{defn}
The \emph{graded Chern character} of $V$ is the characteristic class ${\rm Ch}(V)\in H^\bullet(X)$ defined by ${\rm Ch}(V):=\sum_{j=1}^r\exp(2\pi\sqrt{-1}\delta_j)$.
\end{defn}
\subsection{Statement of the conjecture}Let $X$ be a Fano variety. In \cite{D0} Dubrovin conjectured that many properties of the $qDE$ of $X$, in particular its monodromy, Stokes and central connection matrices, are encoded in the geometry of exceptional collections in $\mathcal D^b(X)$. The following conjecture is a refinement of the original version in \cite{D0}. 
\begin{conj}[\cite{CDG2}]\label{conj}
Let $X$ be a smooth Fano variety of Hodge-Tate type.
\begin{enumerate}
\item The quantum cohomology $QH^\bullet(X)$ has semisimple points if and only if there exists a full exceptional collection in $\mathcal D^b(X)$.
\item If $QH^\bullet(X)$ is generically semisimple, for any oriented ray $\ell$ of slope $\phi\in[0,2\pi[$ there is a correspondence between $\ell$-chambers and helices with a marked foundation.
\item Let $\Omega_\ell$ be an $\ell$-chamber and $\frak E_\ell=(E_1,\dots, E_n)$ the corresponding exceptional collection (the marked foundation). Denote by $S$ and $C$ Stokes and central connection matrices computed in $\Om_\ell$.
\begin{enumerate}
\item The matrix $S$ is the inverse of the Gram matrix of the $\chi$-pairing in $K_0(X)_{\mathbb C}$ wrt the exceptional basis $[\frak E_\ell]$,
\beq
(S^{-1})_{ij}=\chi(E_i,E_j);
\eneq
\item The matrix $C$ coincides with the matrix associated with the $\mathbb C$-linear morphism
\begin{align}
\textnormal{\textcyr{D}}_X^-\colon K_0(X)_\mathbb C\longrightarrow& H^\bullet(X)\\ 
F\xmapsto{\quad\quad}& \frac{(\sqrt{-1})^{\overline{d}}}{(2\pi)^{\frac{d}{2}}}\widehat\Gamma^-_X\exp(-\pi\sqrt{-1}c_1(X)){\rm Ch}(F),
\end{align}
where $d:=\dim_{\mathbb C}X$, and $\overline d$ is the residue class $d\,({\rm mod\,}2)$. The matrix is computed wrt the exceptional basis $[\frak E_\ell]$ and the pre-fixed basis $(T_\alpha)_{\alpha=1}^n$ of $H^\bullet(X)$.
\end{enumerate}
\end{enumerate}
\end{conj}

\begin{rem}
Conjecture \ref{conj} relates two different aspects of the geometry of $X$, namely its \emph{symptectic structure} ($GW$-theory) and its \emph{complex structure} (the derived category $\mathcal D^b(X)$). Heuristically, Conjecture \ref{conj} follows from Homological Mirror Symmetry Conjecture of M.\,Kontsevich, see \cite[Section 5.5]{CDG2}. 
\end{rem}
\begin{rem}
In the paper \cite{KKP} it was underlined the role of $\Gamma$-classes for refining the original version of Dubrovin's conjecture \cite{D0}. Subsequently, in \cite{D3} and \cite[$\Gamma$-conjecture II]{GGI} two equivalent versions of point (3.b) above were given. However, in both these versions, different choices of solutions in Levelt form of the $qDE$ at $z=0$ are chosen wrt the natural ones in the theory of Frobenius manifolds, see Remark \ref{Rc1}, and \cite[Section 5.6]{CDG2}.
\end{rem}
\begin{rem}
If point (3.b) holds true, then automatically also point (3.a) holds true. This follows from the identity \eqref{const2} and Hirzebruch-Riemann-Roch Theorem, see \cite[Corollary 5.8]{CDG2}.
\end{rem}
\begin{rem}
Assume the validity of points (3.a) and (3.b) of Conjecture \ref{conj}. The action of the braid group $\mathcal B_n$ on the Stokes and central connection matrices (Lemma \ref{actbr}) is compatible with the action of $\mathcal B_n$ on the marked foundations attached at each $\ell$-chambers. Different choices of the branch of the $\Psi$-matrix correspond to shifts of objects of the marked foundation. The matrix $M_0^{-1}$ is identified with the canonical operator $\kappa\colon K_0(X)_{\mathbb C}\to K_0(X)_{\mathbb C},\ [F]\mapsto (-1)^{d}[F\otimes \omega_X]$. Equations \eqref{allsc} imply that the connection matrices $C^{(m)}$, with $m\in \mathbb Z$, correspond to the matrices of the morphism \textcyr{D}$^-_X$ wrt the foundations $(\frak E_\ell\otimes \omega_X^{\otimes m})[md]$. The statement  $S^{(m)}=S$ coincides with the periodicity described in Lemma \ref{sper}, see \cite[Theorem 5.9]{CDG2}.
\end{rem}

\begin{rem}
Point (3.b) of Conjecture \ref{conj} allows to identify $K$-classes with solutions of the joint system of equations \eqref{eq1}, \eqref{qde}. Under this identification, Stokes fundamental solutions correspond to exceptional bases of $K$-theory. In the approach of \cite{CV,TV}, where the equivariant case is addressed, such an identification is more fundamental and \emph{a priori}, see Section \ref{qkzsec}.
\end{rem}
\section{Results for Grassmannians}
Conjecture \ref{conj} has been proved for complex Grassmannians $\mathbb G(k,n)$ in \cite{CDG2,GGI}. See also \cite{G1,U}. The proof is based on direct computation of the monodromy data of the $qDE$ at points of the \emph{small quantum cohomology}, namely the subset $H^2(\mathbb G(k,n),\mathbb C)$ of $\Om$. Here we summarize the main results obtained.
\begin{rem}
If\footnote{Here $\pi_1(n)$ denotes the smallest prime number which divides $n$.} $\pi_1(n)\leq k\leq n-\pi_1(n)$,  the small quantum locus of $\mathbb G(k,n)$ is contained in the coalescence locus $\Delta_{\Om}$, see \cite{Cot}. In these cases, the computation of the monodromy data is justified by the results of \cite{CDG0,CDG1,CG1,CG2}. See also Remark \ref{coal}.
\end{rem}
\subsection{The case of projective spaces} Denote by $\sigma\in H^2(\mathbb P^{n-1},\mathbb C)$ the hyperplane class and 
fix the basis $(\sigma^k)_{k=0}^{n-1}$ of $H^\bullet(\mathbb P^{n-1})$. 
The joint system \eqref{eq1}, \eqref{qde} for $\mathbb P^{n-1}$, restricted at the point $t\sigma\in H^2(\mathbb P^{n-1},\mathbb C)$, with $t\in\mathbb C$, is
\begin{align}
\label{eq1p}
\frac{\partial Z}{\partial t}&=z\mathcal C(t)Z,\\
\label{eq2p}
\frac{\partial Z}{\partial z}&=\left(\mathcal U(t)+\frac{1}{z}\mu\right)Z,
\end{align}
with 
{\small
\beq
\mathcal U(t)=\left(\begin{array}{ccccc}
0&&&&nq\\
n&0&&&\\
&n&0&&\\
&& \ddots&\ddots&\\
&&&n&0
\end{array}\right),\quad q:=e^{t},\quad \mathcal C(t)=\frac{1}{n}\mathcal U(t),
\eneq
\beq
\mu=\text{diag}\left(-\frac{n-1}{2},-\frac{n-3}{2},\dots,\frac{n-3}{2},\frac{n-1}{2}\right).
\eneq}
The canonical coordinates are given by the eigenvalues of the matrix $\mathcal U(t)$,
\beq
u_h(t)=ne^{\frac{2\pi i (h-1)}{n}}q^{\frac{1}{n}}\quad h=1,\dots,n.
\eneq
Fix the orthonormalized idempotent vector fields, $f_1(t),\dots, f_n(t)$, given by
\[
f_h(t):=\sum_{\ell=1}^nf_h^\ell(t)\sigma^{\ell-1},\quad f_h^\ell(t):=n^{-\frac{1}{2}}q^{\frac{n+1-2\ell}{2n}}e^{(1-2\ell)i \pi\frac{(h-1)}{n}}\quad h,\ell=1,\dots,n,
\]
and consider the following branch of the $\Psi$-matrix,
\beq\label{psip}
\Psi(t):=\left(\begin{array}{c|c|c}
f_1^1(t)&\dots&f_n^1(t)\\
\vdots&&\vdots\\
f_1^n(t)&\dots& f_n^n(t)
\end{array}\right)^{-1}.
\eneq

\begin{thm}[\cite{CDG2}]\label{conjp}
Fix the oriented ray $\ell$ in $\widehat{\mathbb C^*}$ of slope $\phi\in[0,\frac{\pi}{n}[$. For suitable choices of the signs of the columns of the $\Psi$-matrix \eqref{psip}, the central connection matrix computed at $0\in H^\bullet(\mathbb P^{n-1})$ coincides with the matrix attached to the morphism $$\textnormal{\textcyr{D}}_{\mathbb P^{n-1}}^-\colon K_0(\mathbb P^{n-1})_\mathbb C\to H^\bullet(\mathbb P^{n-1})$$
computed wrt the exceptional bases
{\footnotesize\beq\label{excol1}\mathcal O\left(\frac{n}{2}\right),\bigwedge\nolimits^1\mathcal T\left(\frac{n}{2}-1\right),\mathcal O\left(\frac{n}{2}+1\right),\bigwedge\nolimits^3\mathcal T\left(\frac{n}{2}-2\right),\dots,\mathcal O(n-1),\bigwedge\nolimits^{n-1}\mathcal T
\eneq
}
for $n$ even, and
{\footnotesize
\begin{align}\label{excol2}\mathcal O\left(\frac{n-1}{2}\right),\mathcal O\left(\frac{n+1}{2}\right),&\bigwedge\nolimits^2\mathcal T\left(\frac{n-3}{2}\right),\\
\nonumber
&\mathcal O\left(\frac{n+3}{2}\right),\bigwedge\nolimits^4\mathcal T\left(\frac{n-5}{2}\right),\dots,\mathcal O\left(n-1\right),\bigwedge\nolimits^{n-1}\mathcal T
\end{align}
} for $n$ odd. In particular, Conjecture \ref{conj} holds true for $\mathbb P^{n-1}$. 
\end{thm}

\begin{rem}\label{exccolp}
Exceptional collections \eqref{excol1} and \eqref{excol2} are related to Beilinson's exceptional collection \eqref{bec} by mutations and shifts. For different choices of the ray $\ell$, the exceptional collections attached to the monodromy data computed at $0\in H^\bullet(\mathbb P^{n-1})$ are given (up to shifts) by the following list, see \cite{CDG2,CV}.
\begin{enumerate}
\item {\bf Case $n$ odd:} an exceptional collection either of the form
{\footnotesize
\begin{align*}
\mathcal O\left(-k-\frac{n-1}{2}\right),\ &\mathcal T\left(-k-\frac{n-1}{2}-1\right),\ \mathcal O\left(-k-\frac{n-1}{2}+1\right),\\
\bigwedge\nolimits^3\mathcal T\left(-k-\frac{n-1}{2}-2\right),\ &\mathcal O\left(-k-\frac{n-1}{2}+2\right),
\dots\ ,\ \bigwedge\nolimits^{n-4}\mathcal T\left(-k-n+2\right),\\ \mathcal O(-k-1),\ &\bigwedge\nolimits^{n-2}\mathcal T\left(-k-n+1\right),\ \mathcal O(-k),
\end{align*}
}
or of the form
{\footnotesize
\begin{align*}
\mathcal O\left(-k-\frac{n-1}{2}\right),\ &\mathcal O\left(-k-\frac{n-1}{2}+1\right),\ \bigwedge\nolimits^2\mathcal T\left(-k-\frac{n-1}{2}-1\right),\\ \mathcal O\left(-k-\frac{n-1}{2}+2\right),\ &\bigwedge\nolimits^3\mathcal T\left(-k-\frac{n-1}{2}-2\right)\dots,\ \mathcal O(-k-1),\\
\ \bigwedge\nolimits^{n-3}\mathcal T\left(-k-n+2\right),\ &\mathcal O(-k),\ \bigwedge\nolimits^{n-1}\mathcal T\left(-k-n+1\right),
\end{align*}
}for some $k\in\mathbb Z$
\item {\bf Case $n$ even:} an exceptional collection either of the form
{\footnotesize
\begin{align*}
\mathcal O\left(-k-\frac{n}{2}\right),\ &\mathcal O\left(-k-\frac{n}{2}+1\right),\ \bigwedge\nolimits^2\mathcal T\left(-k-\frac{n}{2}-1\right),\ \mathcal O\left(-k-\frac{n}{2}+2\right),\dots,\\
&\dots\ ,\ \bigwedge\nolimits^{n-4}\mathcal T\left(-k-n+2\right),\ \mathcal O(-k-1),\ \bigwedge\nolimits^{n-2}\mathcal T\left(-k-n+1\right),\ \mathcal O(-k),
\end{align*}
}or of the form 
{\footnotesize
\begin{align*}
\mathcal O\left(-k-\frac{n}{2}+1\right),\ &\mathcal T\left(-k-\frac{n}{2}\right),\ \mathcal O\left(-k-\frac{n}{2}+2\right),\ \bigwedge\nolimits^3\mathcal T\left(-k-\frac{n}{2}-1\right),\dots,\\
\dots&\ ,\ \mathcal O(-k-1),\ \bigwedge\nolimits^{n-3}\mathcal T\left(-k-n+2\right),\ \mathcal O(-k),\ \bigwedge\nolimits^{n-1}\mathcal T\left(-k-n+1\right),
\end{align*}
}for some $k\in\mathbb Z$.
\end{enumerate}
\end{rem}

\subsection{The case of Grassmannians}
Denote by $\mathbb G$ the Grassmannian $\mathbb G(k,n)$ parametrizing $k$-dimensional subspaces in $\mathbb C^n$, and by $\mathbb P$ the projective space $\mathbb P^{n-1}$. Let $\xi_1,\dots, \xi_k$ be the Chern roots of the dual of the tautological bundle $\mathcal S$ on $\mathbb G$, and denote by $h_j(\bm \xi)$ the $j$-th complete symmetric polynomial in $	\xi_1,\dots, \xi_k$. An additive basis of the cohomology ring
\beq\label{prehg}
H^\bullet(\mathbb G)\cong\mathbb C[\xi_1,\dots, \xi_k]^{\frak S_k}\Big/\langle h_{n-k+1},\dots, h_n\rangle,
\eneq
is given by the Schubert classes $(\sigma_\lambda)_{\lambda\subseteq k\times (n-k)}$, labelled by partitions $\lambda$ with Young diagram inside a $k\times(n-k)$ rectangle. Under the presentation \eqref{prehg}, the Schubert classes are given by Schur polynomials in $\bm \xi$,
\beq
\sigma_\lambda:=\frac{\det\left(\xi_i^{\lambda_j+k-j}\right)_{1\leq i,j\leq k}}{\prod_{i<j}(\xi_i-\xi_j)}.
\eneq
Denote by $\eta_{\mathbb P}$ and $\eta_{\mathbb G}$  the Poincar\'e metrics on $H^\bullet(\mathbb P)$ and $H^\bullet(\mathbb G)$ respectively. The metric $\eta_{\mathbb P}$ induces a metric $\eta_{\mathbb P}^{\wedge k}$ on the exterior power $\bigwedge\nolimits^kH^\bullet(\mathbb P)$:
\beq
\eta_{\mathbb P}^{\wedge k}(\alpha_1\wedge\dots,\wedge\alpha_k,\beta_1\wedge\dots,\wedge\beta_k):=\det\left(\eta_{\mathbb P}(\alpha_i,\beta_j)\right)_{1\leq i,j\leq k}.
\eneq
\begin{thm}[\cite{CDG2,GGI}]\label{sat}
We have a $\mathbb C$-linear isometry 
\[
\mathcal I\colon \left(\bigwedge\nolimits^kH^\bullet(\mathbb P),\ (-1)^{\binom{k}{2}}\eta_{\mathbb P}^{\wedge k}\right)\to \left(H^\bullet(\mathbb G),\eta_{\mathbb G}\right),\quad \sigma^{\nu_1}\wedge\dots\wedge\sigma^{\nu_k}\mapsto \sigma_{\tilde\nu},
\]
where $n-1\geq \nu_1>\nu_2>\dots>\nu_k\geq 0$ and $\tilde{\nu}:=(\nu_1-k+1,\nu_2-k+2,\dots,\nu_k)$. 
\end{thm}
Consider the domain $\Om_{\mathbb G}\subset H^\bullet(\mathbb G)$ (resp. $\Om_{\mathbb P}\subset H^\bullet(\mathbb P)$) where the $GW$-potential $F^\mathbb G_0$ (resp. $F^\mathbb P_0$) converges. Let $t\in\mathbb C$ and consider the points 
\beq\label{pp'}
p:=t\sigma_1\in H^2(\mathbb G,\mathbb C),\quad \hat p:=\left(t+\pi\sqrt{-1}(k-1)\right)\sigma\in H^2(\mathbb P,\mathbb C),
\eneq
in the small quantum cohomology of $\mathbb G$ and $\mathbb P$ respectively. Theorem \ref{sat} allow us to identify\footnote{In what follows, if $A$ is a $n\times n$-matrix, we denote by $\bigwedge^kA$ the matrix of $k\times k$-minors of $A$, ordered in lexicographical order.} the tangent spaces $T_p\Om_{\mathbb G}$ and $\bigwedge^kT_{\hat p}\Om_{\mathbb P}$. 

\begin{lem}[\cite{CDG2,GGI}]
Let $\Psi^{\mathbb P}(t)$ be the $\Psi$-matrix defined by \eqref{psip}. Then the matrix $\Psi^{\mathbb G}(t):=(\sqrt{-1})^{\binom{k}{2}}\bigwedge^k\Psi^{\mathbb P}(t+\pi\sqrt{-1}(k-1))$ defines a branch of the $\Psi$-matrix for $\mathbb G$.
\end{lem}

The following results show that under the identification of Theorem \ref{sat}, solutions and monodromy data of the joint system \eqref{eq1}, \eqref{qde} for $\mathbb G$ can be reconstructed from solutions for the joint system for $\mathbb P$.

\begin{thm}[\cite{CDG2}]
Let $Z^{\mathbb P}(t,z)$ be a solution of the joint system \eqref{eq1p}, \eqref{eq2p}. The function
\beq
Z^{\mathbb G}(t,z):=\bigwedge\nolimits^k\left(Z^{\mathbb P}(t+\pi\sqrt{-1}(k-1),z)\right)
\eneq
is a solution for the joint system for $\mathbb G$, namely
\begin{align}
\label{eq1g}
\frac{\partial Z^{\mathbb G}}{\partial t}&=z\mathcal C_{\mathbb G}(t)Z^{\mathbb G},\\
\label{eq2g}
\frac{\partial Z^{\mathbb G}}{\partial z}&=\left(\mathcal U_{\mathbb G}(t)+\frac{1}{z}\mu_{\mathbb G}\right)Z^{\mathbb G}.
\end{align}
\end{thm}

\begin{cor}[\cite{CDG2}]\label{ptog}Fix an oriented ray $\ell$ in $\widehat{\mathbb C^*}$ admissible at both points $p,\hat{p}$ in \eqref{pp'}. Denote by $S^{\mathbb P}(\hat p), S^{\mathbb G}( p)$ and $C^{\mathbb P}(\hat p), C^{\mathbb G}( p)$ the Stokes and central connection matrices at $\hat p$ and $p$, respectively. We have
\begin{align}
S^{\mathbb G}( p)&=\bigwedge\nolimits^kS^{\mathbb P}(\hat p),\\
C^{\mathbb G}( p)&= (\sqrt{-1})^{-\binom{k}{2}}\left(\bigwedge\nolimits^kC^{\mathbb P}(\hat p)\right)\exp(\pi\sqrt{-1}(k-1)\sigma_1\cup).
\end{align}
\end{cor}
\proof
Denote by 
\begin{itemize}
\item $Z^{\mathbb P}_{\rm top}(t,z)$ and $Z^{\mathbb G}_{\rm top}(t,z)$ the topological-enumerative solutions for $\mathbb P$ and $\mathbb G$ respectively, restricted at their small quantum cohomologies;
\item $Z^{\mathbb P/\mathbb G}_{L/R,m}(t,z)$, with $m\in\mathbb Z$, the Stokes fundamental solutions of the joint systems \eqref{eq1}, \eqref{qde} for $\mathbb P$ and $\mathbb G$ respectively.
\end{itemize}
We have
\begin{align*}
Z^{\mathbb G}_{\rm top}(t,z)&=\left(\bigwedge\nolimits^kZ^{\mathbb P}_{\rm top}(t+\pi\sqrt{-1}(k-1),z)\right)\cdot \exp(-\pi\sqrt{-1}(k-1)\sigma_1\cup),\\
Z^{\mathbb G}_{L/R,m}&(t,z)=(\sqrt{-1})^{-\binom{k}{2}}\bigwedge\nolimits^kZ^{\mathbb P}_{L/R,m}(t+\pi\sqrt{-1}(k-1),z).
\end{align*}
See \cite{CDG2} for proofs of these identities.
\endproof

\begin{cor}[\cite{CDG2}]
The central connection matrix computed at $0\in H^\bullet(\mathbb G)$ coincides with the matrix attached to the morphism $$\textnormal{\textcyr{D}}_{\mathbb G}^-\colon K_0(\mathbb G)_\mathbb C\to H^\bullet(\mathbb G)$$
computed wrt an exceptional basis of $K_0(\mathbb G)_\mathbb C$. Such a basis is the projection in $K$-theory of an exceptional collection of $\mathcal D^b(\mathbb G)$ related by mutations and shifts to the twisted Kapranov excptional collection
\beq
(\mathbb S^\lambda\mathcal S^\vee\otimes \mathcal L),\quad \mathcal L:=\det\left(\bigwedge\nolimits^2\mathcal S^\vee\right).
\eneq
In particular, Conjecture \ref{conj} holds true for $\mathbb G$.
\end{cor}

\section{Results on the equivariant $qDE$ of $\mathbb P^{n-1}$}\label{qkzsec}
Gromov-Witten theory, as described in Section \ref{gwsec}, can be suitably adapted to the equivariant case \cite{G}. Given a variety $X$ equipped with the action of a group $G$, a quantum deformation of the equivariant cohomology algebra $H^\bullet_G(X,\mathbb C)$ can be defined. 

Consider the projective space $\mathbb P^{n-1}$ equipped with the diagonal action of the torus $\mathbb T:=(\mathbb C^*)^n$. Although the isomonodromic system \eqref{eq1p}, \eqref{eq2p} does not admit an equivariant analog, the differential equation \eqref{eq1p} only can be easily modified. By change of coordinates $q:=\exp(t)$, setting $z=1$, and replacing the quantum multiplication $*_q$ by the corresponding equivariant one $*_{q,\bm z}$, equation \eqref{eq1p} takes the form
\beq\label{qde.2}
q\frac{d}{dq}Z=\sigma*_{q,\bm z}Z.
\eneq
Here the equivariant parameters $\bm z=(z_1,\dots, z_n)$ correspond to the factors of $\mathbb T$, and $Z(q,\bm z)$ takes values in $H^\bullet_{\mathbb T}(\mathbb P^{n-1},\mathbb C)$. Equation \eqref{qde.2} admits a compatible system of difference equations, called $qKZ$ difference equations
\beq\label{qkz}
Z(q,z_1,\dots,z_i-1,\dots,z_n)=K_i(q,\bm z)Z(q,\bm z),\quad i=1,\dots,n,
\eneq
for suitable linear operators $K_i$'s, introduced in \cite{TV}. The joint system \eqref{qde.2}, \eqref{qkz} is a suitable limit of an analogue one for the cotangent bundle $T^*\mathbb P^{n-1}$, see \cite{GRTV,RTV}. The existence and compatibility of such a joint system for more general Nakajima quiver varieties is justified by the general theory of D.\,Maulik and A.\,Okounkov \cite{MO}.

In \cite{TV}, the study of the monodromy and Stokes phenomenon at $q=\infty$ of solutions of the joint system \eqref{qde.2}, \eqref{qkz} is addressed. Furthermore, elements of $K_0^{\mathbb T}(\mathbb P^{n-1})_{\mathbb C}$ are identified with solutions of the joint system \eqref{qde.2}, \eqref{qkz}: Stokes bases of solutions correspond to exceptional bases.

In \cite{CV}, the authors describe relations between the monodromy data of the joint system of the equivariant $qDE$ \eqref{qde.2} and $qKZ$ equations \eqref{qkz} and characteristic classes of objects of the derived category $\mathcal D^b_{\mathbb T}(\mathbb P^{n-1})$ of equivariant coherent sheaves on $\mathbb P^{n-1}$. Equivariant analogs of results of \cite[Section 6]{CDG2} are obtained. 

The $\textcyr{B}$-Theorem of \cite{CV} is the equivariant analog of Theorem \ref{conjp}. Moreover, in \cite{CV} the Stokes bases of solutions of the joint system \eqref{qde.2}, \eqref{qkz} are identified with explicit $\mathbb T$-full exceptional collections in $\mathcal D^b_{\mathbb T}(\mathbb P^{n-1})$, which project to those listed in Remark \ref{exccolp} via the forgetful functor $\mathcal D^b_{\mathbb T}(\mathbb P^{n-1})\to\mathcal D^b(\mathbb P^{n-1})$. This refines results of \cite{TV}. Finally, in \cite{CV} it is proved that the Stokes matrices of the joint system \eqref{qde.2}, \eqref{qkz} equal the Gram matrices of the equivariant Grothendieck-Euler-Poincar\'e pairing on $K_0^\mathbb T(\mathbb P^{n-1})_{\mathbb C}$ wrt the same exceptional bases.

% ------------------------------------------------------------------------

\subsection*{Acknowledgment}
These notes partly touch the topic of the talk given by the author at the XXXVIII \emph{Workshop on Geometric Methods in Physics}, hold in June-July 2019 in the inspiring atmosphere of Bia\l{}owie\.{z}a, Poland. The author is thankful to the organizers of the Workshop for invitation. He also thanks the Max-Planck-Institut f\"ur Mathematik in Bonn, Germany, for support.

% ------------------------------------------------------------------------
\end{document}